\documentclass{elsart}
\usepackage{amsfonts}
\usepackage{fancyhdr}
\usepackage[latin1]{inputenc}

\usepackage{natbib}
\usepackage{fancyhdr}
\usepackage{amsmath}
\usepackage{graphicx}
\usepackage{hyperref}
\usepackage{amssymb}
\usepackage{dsfont}
\newtheorem{theorem}{Theorem}[section]
\newtheorem{remark}{Remark}[section]
\newtheorem{lemma}{Lemma}[section]
\numberwithin{equation}{section}

\renewcommand{\cite}{\citet}

\usepackage{amssymb}

\begin{document}

\begin{frontmatter}



\title{Testing additivity in nonparametric regression under random censorship}

\author[k]{Mohammed DEBBARH}
\address[k]{L.S.T.A.-Universit\'{e} Paris 6, 175, rue du Chevaleret,
 75013, Paris. debbarh@ccr.jussieu.fr}

\author[k,l]{Vivian  VIALLON}
\address[l]{Department of Biostatistics, H\^opital Cochin, Universit\'e Paris-Descartes, 27 rue du Faubourg Saint Jacques 75014 Paris. vivian.viallon@univ-paris5.fr}
\begin{abstract}
In this paper, we are concerned with nonparametric estimation of
the multivariate regression function in the presence of right
censored data. More precisely, we propose a statistic that is
shown to be asymptotically normally distributed under the additive
assumption, and that could be used to test for additivity in the
censored regression setting.
\end{abstract}

\begin{keyword}
censored regression; additive model; curse of dimensionality;
nonparametric regression; marginal integration.
\end{keyword}
\end{frontmatter}

\vspace{-1cm}
\section{Introduction and motivations}
\vspace{-0.6cm} A well known issue in nonparametric regression
estimation is the so-called \emph{curse of dimensionality}, i.e.
the fact that the rate of convergence of nonparametric estimators
dramatically decreases as the dimension of the covariates
increases (see, for instance, \cite{Stone1982}). To get round this
issue, one common solution is to work under the \emph{additive
assumption}, i.e. the true regression function is assumed to be
the sum of some lower dimension regression functions (typically,
univariate or bivariate functions). But, this assumption is strong
and has therefore to be checked \emph{via} one of the available
tests (\cite{Camlong2001}, \cite{Linton2001}, \cite{Sperlich2002},
\cite{Derbort}) before being used in practice.

When the variable of interest is censored, several nonparametric
estimators have been proposed for the multivariate regression
function (see, e.g., \cite{FanGijbels}, \cite{Carbonez1995},
\cite{Kohler2002}, \cite{BrunelComte1}). By combining one of this
'initial' estimator with the marginal integration method (see
\cite{newey}, \cite{Linton1995}), estimates can be obtained under
the additive assumption. In particular, \cite{DebbarhViallon} made
use of an initial \emph{Inverse Probability of Censoring Weighted}
estimator (such as the one proposed by \cite{Carbonez1995}), and
established the uniform convergence rate for the corresponding
additive estimator. However, in this censored setting, no test for
additivity has been proposed yet. That will be our concern here.
Namely, we first exhibit a statistic evaluating a weighted
difference between the observations of the variable of interest
and the estimator we derive \emph{via} the marginal integration
method. Then, this statistic is shown to be asymptotically
normally distributed under the additive assumption.

To build our estimators, and then our test statistic, the
following notations are needed. Let $(Y, C, {\bf X}),$ $(Y_1, C_1,
{\bf X}_1),$ $(Y_2, C_2, {\bf X}_2),...$ be independent and
identically distributed $\mathbb{R} \times \mathbb{R} \times
\mathbb{R}^d$-valued random variables. Here $Y$ is the variable of
interest, $C$ the censoring variable and ${\bf X}=(X_1,...,X_d)$ a
vector of concomitant variables. In the right censorship model,
the only available information on $(Y,C)$ is given by
$(Z,\delta)$, with $Z=\min\{Y,C\}$ and $\delta=\mathbb{I}_{\{{Y}
\leq C\}}$, $\mathbb{I}_E$ standing for the indicator function of
the set $E$. As a matter of fact, the observed sample is
${\mathcal D}_n:=({\bf X}_i,Z_i,\delta_i)_{1\leq i \leq n}$, for a
given $n\geq 1$. \vskip5pt \noindent Given a real measurable
function $\psi$, our concern here is the regression function of
$\psi(Y)$ evaluated at $\bf{X}=\bf{x}$, that is,
\vspace{-.5cm}
\begin{eqnarray}\label{add_component}
m_{\psi}(\bf{x}) &=&  E\left(\psi (Y) \mid \bf{X}=\bf{x}\right),~
\forall ~{\bf x}=(x_1,...,x_d) \in \mathbb{R}^d. \label{fdereg}
\end{eqnarray}
Under the traditional additive assumption, the regression function
defined in (\ref{fdereg}) can be written as the sum of some
(unknown) univariate regression functions $m_l$, \vspace{-.5cm}
\begin{eqnarray}
m_{\psi}({\bf x})=m_{\psi,add}({\bf x}):=\mu +
\sum_{l=1}^d m_l(x_l)\label{add_model}.
\end{eqnarray}
In view of (\ref{add_model}), the functions $m_l$, as well as the
constant term $\mu$, are defined up to an additive constant.
Therefore, we will work under the common identifiability condition
$Em_l(X_l)=0$, for $l=1,...,d$. This condition implies that
$\mu=E(\psi(Y))$. \vskip5pt
 \noindent In the sequel, we set, for
all $t\in \mathbb{R}$, $F(t)=P(Y> t)$, $G(t)=P(C> t)$ and
$H(t)=P(Z>t)$ the survival functions pertaining to $Y$, $C$ and
$Z$ respectively. Further denote by $G_n$ the Kaplan-Meier
(\cite{KM}) estimator of $G$, \space{-.5cm}
\begin{equation}
G_n (y) =\prod_{1\leq i \leq
n}\!\!\!{\Big(\frac{N_n(Z_i)-1}{N_n(Z_i)}\Big),}^{\!\!\!\!
\beta_i}\mbox{ for all }y\geq 0,\mbox{ with }
\beta_i=\mathbb{I}_{\{Z_i\leq y\}}(1-\delta_i). \label{Gnstar}
\end{equation}
Here we defined $N_n(x)=\sum_{i=1}^n\mathbb{I}_{\{Z_i\leq x\}}$,
and the conventions $\prod_{\varnothing}=1$ and $0^0=1$ were
adopted. \vskip5pt \noindent Consider the null hypothesis
\vspace{-0.5cm}
\begin{equation*}
H_0:~~ m_\psi\in \mathcal{M}_{add}:=\{m:\mathbb{R}^d\rightarrow \mathbb{R}, m({\bf
x})=\mu +\sum_{l=1}^dm_l(x_l); E(m_l(X_l))=0\}.
\end{equation*}
Following the ideas of \cite{Mammen}, \cite{Camlong2001} and
\cite{Vieu}, we denote by $g$ some fixed weight function, by $L$ a
given \emph{kernel}, i.e. a real measurable function integrating
to 1, defined in $\mathbb{R}^d$ and by $(\ell_n)_{n\geq1}$ a
sequence of positive real numbers. Further let $\widehat
m_{\psi,add}^\star$ be some estimator of $m_\psi$ under the
additive assumption (\ref{add_model}) (or, equivalently, under
$H_0$). Now, let us consider the statistic
\vspace{-0.5cm}
\begin{eqnarray}\label{est_stat_test}
T_n^\star = \int_{\mathbb{R}^d}\Big[
\frac{1}{n\ell_n^d}\sum_{i=1}^n L\Big(\frac{{\bf x}-{\bf
X}_i}{\ell_n}\Big)\Big(\frac{\frac{\delta_i\psi(Z_i)}{G_n(Z_i)}-
\widehat m_{\psi,add}^\star({\bf X}_i)}{\hat f_n({\bf X}_i)}\Big)
\Big]^2g({\bf x})d{\bf x},
\end{eqnarray}
which is a natural estimator of the quantity
$E(E(\frac{\delta\psi(Z)}{G(Z)}-m_{\psi,add}({\bf X}))|{\bf
X})^2$. Under a useful independence condition (see $(C.1)$ below),
the latter quantity equals $E(E(\psi(Y)-m_{\psi,add}({\bf
X}))|{\bf X})^2$, and then equals zero if and only if the
hypothesis $H_0$ is true. Moreover, in Theorem \ref{th_norm_stat}
below, this statistic is shown to be asymptotically normally
distributed under $H_0$. Therefore, it could be useful to test for
additivity in censored nonparametric regression. Properties of the
corresponding test will be studied elsewhere.\vskip5pt \noindent
Now, we precise how $\widehat m_{\psi,add}^\star$ may be
constructed. Let $K_1$, $K_2$, $K_3$ and $K$, be kernels
respectively defined in $\mathbb{R}$, $\mathbb{R}^{d-1}$,
$\mathbb{R}^d$ and $\mathbb{R}^d$. Further set $\hat f_n$ the
kernel estimator of $f$, with $f$ denoting the density function of
${\bf X }$. Namely, \vspace{-0.5cm}
\begin{eqnarray*}
\hat f_n({\bf x})= \frac{1}{nh_n^d} \sum_{j=1}^n K\Big(\frac{{\bf
X}_j-{\bf x}}{h_n}\Big),
\end{eqnarray*}
where $(h_n)_{n \geq 1}$ is a given sequence of positive real
numbers. Denote by $(h_{j,n})_{n \geq 1}$, $j=1,2$, two sequences
of positive real numbers. To estimate the multivariate regression
function defined in (\ref{fdereg}), the following Nadaraya-Watson
type estimators can be used (see \cite{Carbonez1995},
\cite{Kohler2002} and \cite{Jones1994}), \vspace{-0.5cm}
\begin{equation} \widetilde{m}_{\psi,n}^\star
({\bf x}) = \sum_{i=1}^n W_{n,i}({\bf x})\frac{\delta_i
\psi(Z_i)}{G_n(Z_i)} ~~ \mbox{with}~~W_{n,i}({\bf x})= \frac{
K_{3}\big({\frac{{\bf x} -{\bf
X}_{i}}{h_{1,n}}}\big)}{nh_{1,n}^d\hat f_n ({\bf
X}_i)},\label{estcarbo1}
\end{equation}
and, for $l=1,...,d,$
\begin{equation}
\widetilde{m}_{\psi,n,l}^\star ({\bf x}) = \sum_{i=1}^n
W_{n,i}^l({\bf x})\frac{\delta_i \psi(Z_i)}{G_n(Z_i)} ~~
\mbox{with}~~W_{n,i}^l({\bf
x})=\frac{K_1\big(\frac{x_l-X_{i,l}}{h_{1,n}}\big)
K_2\big(\frac{{\bf x}_{-l}-{\bf
X}_{i,-l}}{h_{2,n}}\big)}{nh_{1,n}h_{2,n}^{d-1} \hat f_n ({\bf
X}_{i})},\label{estcarbo}
\end{equation}
where we set, for all ${\bf x}=(x_1,..,x_{d})\in\mathbb{R}^d$ and
every $l=1,...,d$, ${\bf x}_{-l}=(x_1,..,$ $x_{l-1},x_{l+1},$
$..,$ $x_d)$.To estimate the additive components, we use the
marginal integration method (see \cite{newey} or
\cite{Linton1995}). Let $q_1,...,q_d$ be $d$ given density
functions. Then, setting $q({\bf x}) = \prod_{l=1}^d q_l(x_l)$ and
$q_{-l}({\bf x}_{-l}) =$ $\prod_{j \neq l} q_j(x_j)$, we define
\vspace{-0.5cm}
\begin{eqnarray}
\eta_{l}(x_{l}) = \int_{\mathbb{R}^{d-1}} m_{\psi}({\bf x})
q_{-l}({\bf x}_{-l}) d{\bf x}_{-l} - \int_{\mathbb{R}^d}
m_{\psi}({\bf x}) q({\bf x}) d{\bf x},\quad l=1,...,d,
\label{additive_component}
\end{eqnarray}
in such a way that the two following equalities hold,
\vspace{-0.5cm}
\begin{eqnarray}
& &\eta_l(x_l) = m_l(x_l) - \int_{\mathbb{R}} m_l(z)
q_l(z)dz, \quad l=1,...,d, \label{relation_additive_component}\\
&&m_{\psi}({\bf x})= \sum_{l=1}^d \eta_l(x_l) +
\int_{\mathbb{R}^d} m_{\psi}({\bf z})q({\bf z}) d{\bf z}
\label{additive_component_marginale}.
\end{eqnarray}
In view of (\ref{relation_additive_component}) and
(\ref{additive_component_marginale}), the functions $\eta_l$,
$l=1,...,d,$ turn out to be some additive components, and, from
(\ref{estcarbo}) and (\ref{additive_component}), a natural
estimator of the $l$-th component $\eta_l$ is given, for all
$l=1,...,d$, by \vspace{-0.5cm}
\begin{eqnarray}
\widehat \eta^\star_{l}(x_{l}) = \int_{\mathbb{R}^{d-1}}\!
\widetilde{m}_{\psi,n,l}^{\star}({\bf x}) q_{-l}({\bf x}_{-l})
d{\bf x}_{-l} - \int_{\mathbb{R}^d}\!
\widetilde{m}_{\psi,n,l}^{\star}({\bf x}) q({\bf x}) d{\bf x}.
\label{est_comp_add}
\end{eqnarray}
From (\ref{est_comp_add}), an estimator $\widehat
m_{\psi,add}^\star$ of the censored regression function can be
deduced under the additive assumption (\ref{add_model}) (or,
equivalently, $H_0$), \vspace{-0.5cm}
\begin{eqnarray}
\widehat m_{\psi,add}^\star({\bf x}) & = & \sum_{l=1}^d \hat
\eta^\star_l(x_l) + \int_{\mathbb{R}^d}
\widetilde{m}_{\psi,n}^{\star}({\bf x}) q({\bf x})d{\bf x}
\label{estim_add}.
\end{eqnarray}
\section{Hypotheses and Results}
\vspace{-0.5cm} \noindent These preliminaries being given, we
introduce the assumptions to be made to state our results. First,
consider the hypotheses pertaining to $(Y,C,{\bf X})$. We suppose
that $({\bf X},Y)$ has a joint density $f_{{\bf X},Y}$. Moreover,
we impose the following conditions. \vspace{-0.5cm}
\begin{eqnarray*}
(C.1):&  & C ~~\mbox{and}~~ ({\bf X},Y) ~~\mbox{are independent}. \\
(C.2):&  &\mbox{G is continuous}.\\
(C.3):&  & \mbox{There exists a constant}~M<\infty~~\mbox{such
that
$\sup_{0\leq t\leq \tau}|\psi(t)| \leq M$}.\\
(C.4):&  & m_\psi \mbox{ is a $k$-times continuously
differentiable function, $ k \geq 1$, and}\\&&~ \tiny{\sup_{\bf
x}\Big|\frac{\partial^km_{\psi}}{\partial x_l^k} ({\bf x})
\Big|<\infty};~l=1,...,d.
\end{eqnarray*}
\vspace{-0.5cm}
\begin{rem}
It is noteworthy that condition $(C.1)$
is stronger than the conditional independence of $C$ and $Y$ given
${\bf X}$, under which \cite{Beran} worked to build an estimator of
the conditional survival function (see also \cite{Dabrowska}).
Note, however, that the two assumptions coincide if
$C$ and ${\bf X}$ are independent. In other respect, to use
Beran's local Kaplan-Meier estimator, the censoring has to be
locally fair, that is $P[C\geq t\mid {\bf X}={\bf x}]>0$
whenever $P[Y\geq t\mid {\bf X}={\bf x}]>0$. Here (see assumption $({\bf A})(ii)(b)$ below), we
essentially suppose that $G(t)>0$ whenever $F(t)>0$, which is,
on its turn, a weaker assumption. For a nice discussion on the
difference between Beran's estimator and Inverse Probability of Censoring Weighted type estimators,
we refer to \cite{Carbonez1995}.
\end{rem}
\noindent Denote by $\mathcal{C}_1, ...,~ \mathcal{C}_d$, $d$
compact intervals of $\mathbb{R}$ and set
$\mathcal{C}=\mathcal{C}_1\times...\times \mathcal{C}_d$. For
every subset $\mathcal{E}$ of $\mathbb{R}^q$, $q\geq1$, and any
$\alpha>0$, introduce the $\alpha$-neighborhood
$\mathcal{E}^\alpha$ of $\mathcal{E}$, i.e. $\mathcal{E}^\alpha =
\{x : \inf_{y\in \mathcal{E}}\|x-y\|_{\mathbb{R}^q}\leq \alpha\}$,
$\|\cdot\|_{\mathbb{R}^q}$ standing for the euclidian norm on
$\mathbb{R}^q$.\\ We will work under the following regularity
assumptions on $f$ and $f_l$, $l=1, ..., d$, $f_l$ denoting the
density function of $X_l$. These functions are supposed to be
continuous and we assume the existence of a constant $\alpha>0$
such that the following assumptions hold, \vspace{-0.5cm}
\begin{eqnarray*}
(F.1):&  & \forall x_l\in\mathcal{C}^\alpha_l, f_{l}(x_l)>0,\
l=1, ..., d, \mbox{ and } \forall{\bf x}\in\mathcal{C}^\alpha f({\bf x})>0.\  \\
(F.2):&  & f \mbox{ is $k'$-times continuously differentiable on }
\mathcal{C}^\alpha, k'> kd.
\end{eqnarray*}
Regarding the weight function $g$, we will assume that the
condition $(G.1)$ below is satisfied. \vspace{-0.5cm}
\begin{eqnarray*}
(G.1):& & g \mbox{ is an indicator function with compact support
included in } \mathcal{C}.
\end{eqnarray*}
\vspace{-0.5cm}
\begin{remark}
Assumption $(G.1)$ is made here to avoid technical issues in the
derivation of our results. Moreover, it is not restrictive since
$g$ is a given weight function. That being said, this assumption
could be relaxed (see, e.g., \cite{Vieu}, \cite{Camlong2001}).
\end{remark}
The kernels $L$, $K$, and $K_3$ defined in $\mathbb{R}^d$, $K_1$
defined in $\mathbb{R}$ and $K_2$ defined in $\mathbb{R}^{d-1}$,
are assumed to be continuous, compactly supported and integrating
to 1. Moreover, we suppose that, \vspace{-0.5cm}
\begin{eqnarray*}
(K.1):&  & K_1 \mbox{ is Lipschitz};\\
(K.2):&  &  \mbox{$K_1$ and $K_3$ are of order $k$, and $K$ is of
order $k'$.}
\end{eqnarray*}
In addition, we impose the following assumptions on the
integrating density functions $q_{-l}$ and $q_l$, $l=1, ..., d$.
\vspace{-0.5cm}
\begin{eqnarray*}
(Q.1):&  & q_{-l} \mbox{ is bounded and
continuous, } l=1, ..., d.\\
(Q.2):&  & q_l \mbox{ has } k+1 \mbox{ continuous and bounded
derivatives, }l=1, ..., d.
\end{eqnarray*}
Turning our attention to the smoothing parameters $\ell_n$, $h_n$
and $h_{j,n},~j=1,2$, we will work under the conditions below.
\vspace{-0.5cm}
\begin{eqnarray*}
(H.1):&  & h_n=c_1\Big(\frac{\log n}{n}\Big)^{1/(2k'+d)}, \mbox{
for a given } 0<c_1<\infty.\\
(H.2):&  & h_{1,n}= c_2\Big(\frac{\log n}{n}\Big)^{1/(2k+1)},
\mbox{ for a given } 0<c_2<\infty  \mbox{ and }~h_{2,n}=o(1).\\
(H.3):&  & n (\log n /n)^{k/(2k+1)}\ell_n^{d/2}\rightarrow 0
~\mbox{and}~n\ell_n^{d} \rightarrow \infty.
\end{eqnarray*}
\noindent As mentioned in \cite{GrossLai}, functionals of the
(conditional) law can generally not be estimated on the complete
support when the variable of interest is right-censored.
Accordingly, we will work under the assumption $({\bf A})$ that
will be said to hold if either $({\bf A})(i)$ or $({\bf A})(ii)$
below holds. Denote by $T_L=\sup\{t:L(t)>0\}$ the upper endpoint
of the distribution of a random variable with right continuous
survival function $L$.
\begin{tabbing}
$({\bf A})(i)\;\;$ \=\ There exists
a $\tau_0<T_H$ such that  $\psi=0$ on $(\tau_0,\infty)$.\\[0.2cm]
$({\bf A})(ii)$ \>\ $(a)\;\;$ For a given $k/(2k+1)<p\leq 1/2$,
$\big|\int_0^{T_H}
 F^{-p/(1-p)}dG\big|<\infty$;\\
\>\ $(b)\;\;$ $T_F<T_G$;\\
\>\ $(c)\;\;$ $ n^{2p-1}h^{-1}_{l,n}|\log(h_{l,n})|
\rightarrow \infty$, as $n\rightarrow\infty$, for every
$l=1,...,d$.
\end{tabbing}
It is noteworthy that assumption $({\bf A})(ii)$ allows for
considering the estimation of the "classical" regression function,
which corresponds to the choice $\psi(y)=y$. On the other hand,
normality for estimators of functionals such as the conditional
distribution function $P(Y\leq \tau_0|{\bf X})$ can be obtained
under weaker conditions, when restricting ourselves to
$\tau_0<T_H$.\vskip5pt
\noindent To state our result, some
additional notations are needed. Set
$\epsilon_i=\frac{\delta_i\psi(Z_i)} {G(Z_i)}-m_\psi({\bf X}_i)$
and $\sigma_0^2({\bf x}) = E(\epsilon_i^2|{\bf X}_i={\bf x})$.
Further introduce $B=[\int \sigma_0^2({\bf u})f^{-1}({\bf
u})g({\bf u})$ $d{\bf u}]\times \int L^2({\bf t})d{\bf t}$ and
$V=2\left[\int (\sigma_0^2({\bf u})^2 f^{-2}({\bf u}) g^2({\bf u})
d{\bf u}\right]\times \int\left[\int L({\bf t})L({\bf t}-{\bf
r})d{\bf t}\right]^2d{\bf r}$.
\begin{theorem}\label{th_norm_stat}
Assume the conditions $({\bf A})$, $(C.1$-$4)$, $(F.1$-$2)$,
$(G.1)$, $(K.1$-$2)$, $(Q.1$-$2)$ and $(H.1$-$3)$ hold. Then,
under the null hypothesis $H_0$, we have,
\begin{eqnarray*}
\frac{n\ell_n^{d/2}T_n^\star -
B\ell_n^{-d/2}}{\sqrt{V}}\rightarrow \mathcal{N}(0,1), \quad\mbox{
as }\ n\rightarrow\infty.
\end{eqnarray*}
\end{theorem}
\vspace{-0.5cm}
\section{Proof of Theorem \ref{th_norm_stat}}
\vspace{-0.5cm} Here, we present the proof of Theorem
\ref{th_norm_stat} in the case where $({\bf A})(i)$ holds. The
case where $({\bf A})(ii)$ holds follows from similar arguments
(especially replacing the result of \cite{Folder1981} by that of
\cite{GuLai} or that of \cite{ChenLo}); details are then omitted.

We will make frequent use of the following lemma, which was established
in \cite{DebbarhViallon}.
\begin{lemma}\label{lem_unif_consist}
Assume $H_0$ holds. Then, under the conditions $({\bf A})$, $(C.1$-$4)$,
$(F.1$-$2)$, $(K.1$-$2)$, $(Q.1$-$2)$ and $(H.1$-$2)$, we have,
with probability one,
\begin{eqnarray}
\sup_{{\bf x} \in \mathcal{C}} |\widehat{m}_{\psi,add}^\star ({\bf
x}) - m_{\psi}({\bf x})| = \mathcal{O}\Big(\Big(\frac{\log
n}{n}\Big)^{\frac{k}{2k+1}}\Big). \label{res_lem_unif_consist}
\end{eqnarray}
\end{lemma}

We will also make frequent use of the following result, due to
\cite{Folder1981}. \vspace{-0.5cm}
\begin{eqnarray}\label{LILfoldes}
\mbox{For all}\,\ \tau'<T_H,\quad \sup_{y\leq \tau'}
|G^\star_n(y)-G(y)|=\mathcal{O}((\log \log n/n)^{1/2})=:\rho_n
\end{eqnarray}
From this last result, we especially get the following type of
approximations. Set $\epsilon_i^\star
=\frac{\delta_i\psi(Z_i)}{G_n(Z_i)} - m_{\psi}({\bf X}_i)$.  Then,
from $({\bf A})(i)$, $(C.2)$, $(C.3)$ and (\ref{LILfoldes}), we
have, almost surely as $n\rightarrow\infty$, \vspace{-0.5cm}
\begin{eqnarray}
\epsilon_i^\star = \epsilon_i + \mathcal{O}(\rho_n).
\label{argument2}
\end{eqnarray}
Now, recalling the definition (\ref{est_stat_test}) of
$T_n^\star$, we have \vspace{-0.5cm}
\begin{eqnarray*}
T_n^\star = \int_{\mathbb{R}^d}\Big[\frac{1}{n\ell_n^d}
\sum_{i=1}^n L\big(\frac{{\bf x}- {\bf
X}_i}{h_n}\Big)\Big(\frac{m_{\psi}({\bf X}_i)- \widehat
m_{\psi,add}^\star({\bf X}_i)+\epsilon_i^\star}{f_n({\bf
X}_i)}\Big)\Big]^2g({\bf x}) d{\bf x}.
\end{eqnarray*}
Consider the quantity \vspace{-0.5cm}
\begin{eqnarray*}
T_n^{1^\star}&=&\int_{\mathbb{R}^d}\Big[\frac{1}{n\ell_n^d}
\sum_{i=1}^n L\Big(\frac{{\bf x}- {\bf
X}_i}{\ell_n}\Big)\Big(\frac{m_{\psi}({\bf X}_i)- \widehat
m_{\psi,add}^\star({\bf X}_i)+\epsilon_i^\star}{f({\bf
X}_i)}\Big)\Big]^2g({\bf
x}) d{\bf x}. \\
&=& T_{n,1}^{1^\star}+ T_{n,2}^{1^\star}+ T_{n,3}^{1^\star}+
2T_{n,4}^{1^\star}+2 T_{n,5}^{1^\star},
\end{eqnarray*}
with (see \cite{Camlong2001}), \vspace{-0.5cm}
\begin{eqnarray*}
T_{n,1}^{1^\star}&= &\int_{\mathbb{R}^d}
\frac{1}{n^2\ell_n^{2d}}\sum_{i=1}^n L^2\Big(\frac{{\bf x}- {\bf
X}_i}{\ell_n}\Big)\frac{{\epsilon_i^\star}^2
g({\bf x})}{f^2({\bf X}_i)}d{\bf x},\\
T_{n,2}^{1^\star}&= &\int_{\mathbb{R}^d}
\frac{1}{n^2\ell_n^{2d}}\sum_{i\neq j} L\Big(\frac{{\bf x}- {\bf
X}_i}{\ell_n}\Big)L\Big(\frac{{\bf x}- {\bf
X}_j}{\ell_n}\Big)\frac{\epsilon_i^\star \epsilon_j^\star g({\bf
x})}{f({\bf X}_i)f({\bf
X}_j)}d{\bf x},\\
T_{n,3}^{1^\star}&= &\int_{\mathbb{R}^d}
\left[\frac{1}{n\ell_n^{d}}\sum_{i=1}^n L\Big(\frac{{\bf x}- {\bf
X}_i}{\ell_n}\Big)\Big(\frac{m_{\psi}({\bf X}_i)-\widehat
m_{\psi,add}^\star({\bf
X}_i)}{f({\bf X}_i)}\Big)\right]^2 g({\bf x})d{\bf x},\\
T_{n,4}^{1^\star}&= &\int_{\mathbb{R}^d}
\left[\frac{1}{n\ell_n^{d}}\right]^2\sum_{i=1}^n
L^2\Big(\frac{{\bf x}- {\bf X}_i}{\ell_n}\Big)(m_{\psi}({\bf
X}_i)-\widehat m_{\psi,add}^\star({\bf
X}_i))\epsilon_i^\star \frac{g({\bf x})}{f^2({\bf X}_i)}d{\bf x},\\
T_{n,5}^{1^\star}&= &\int_{\mathbb{R}^d}
\left[\frac{1}{n\ell_n^{d}}\right]^2\sum_{i\neq j}
L\Big(\frac{{\bf x}- {\bf X}_i}{\ell_n}\Big) L\Big(\frac{{\bf x}-
{\bf X}_j}{\ell_n}\Big)\Big(\frac{m_{\psi}({\bf X}_i)-\widehat
m_{\psi,add}^\star({\bf X}_i)}{f({\bf X}_i)f({\bf
X}_j)}\Big)\epsilon_j^\star g({\bf x})d{\bf x}.
\end{eqnarray*}
By $(F.1$-$2)$, $(H.1)$ and $(K.2)$, it holds that, almost surely
as $n\rightarrow\infty$ (see for instance \cite{Ango-Nze}),
\vspace{-0.5cm}
\begin{eqnarray}
\sup_{{\bf x}\in \mathcal{C}} \big|\hat f_n ({\bf x})-f({\bf
x})\big|=\mathcal{O}\Big(\sqrt{\frac{\log n}{nh_n^d}}\Big),
\label{ango}
\end{eqnarray}
in such a way that  $T_n^\star = T_n^{1^\star}\{1+
\mathcal{O}[(n^{-1}h_n^{-d}\log n)^{1/2}]\}$ almost surely as
$n\rightarrow\infty$. Therefore, to achieve the proof of Theorem
\ref{th_norm_stat}, it is sufficient to establish
(\ref{1})--(\ref{5}) below. \vspace{-0.5cm}
\begin{eqnarray}
& &n\ell_n^dT_{n,1}^{1^\star}=B+o_p(\ell_n^{d/2}),\label{1}
\\ & &n\ell_n^{d/2} \frac{(T_{n,2}^{1^\star}-ET_{n,2}^{1^\star})}{\sqrt{V}} \rightarrow
N(0,1),\label{2}\\
&  & T_{n,3}^{1^\star}=o_p(n^{-1}\ell_n^{-d/2}),\label{3}\\
&  &T_{n,4}^{1^\star}=o_p(n^{-1}\ell_n^{-d/2}),\label{4}\\
& &T_{n,5}^{1^\star}=o_p(n^{-1}\ell_n^{-d/2}).\label{5}
\end{eqnarray}
{\bf Proof of (\ref{1}):}~ Set ${\sigma_0^{\star}}^2({\bf x}) =
E(\epsilon^{\star\ 2}_i|{\bf X}_i={\bf x})$. Using a conditioning
argument, it is straightforward that \vspace{-0.5cm}
\begin{eqnarray*}
ET_{n,1}^{1^\star}&= & \frac{1}{n\ell_n^{2d}}\int \int
\frac{{\sigma_0^\star}^2({\bf v})}{f({\bf v})}L^2\Big(\frac{{\bf
x}- {\bf v}}{\ell_n}\Big) d{\bf v} g({\bf x}) d{\bf x}.
\end{eqnarray*}
Moreover, since $g$ is an indicator function with compact support
included in ${\mathcal C}$, we obtain, for $n$ large enough, that
\vspace{-0.5cm}
\begin{eqnarray*}
ET_{n,1}^{1^\star}& = &\frac{1}{n\ell_n^d}\int
E({\epsilon_1^\star}^2|{\bf X}_1={\bf x})g({\bf x})f^{-1}({\bf
x})d{\bf x} \int L^2({\bf u})d{\bf
u}+\mathcal{O}\Big(\frac{1}{n}\Big),
\end{eqnarray*}
and then, \emph{via} arguments similar to those used to derive
(\ref{argument2}), \vspace{-0.5cm}
\begin{eqnarray*}
ET_{n,1}^{1^\star}& =& \frac{1}{n\ell_n^d}B+
\mathcal{O}(\rho_n),\quad \mbox{ as $n\rightarrow\infty$}.
\end{eqnarray*}
Turning our attention to the variance of $T_{n,1}^{1^\star}$, we
can write, \vspace{-0.5cm}
\begin{eqnarray}\label{vartn1}
{\rm Var }\ T_{n,1}^{1^\star}= \frac{1}{n^3\ell_n^{4d}}\ {\rm Var
}(I), \mbox{ where } I=\int
{\epsilon_1^\star}^2L^2\left(\frac{{\bf x}-{\bf
X}_j}{\ell_n}\right)\frac{g({\bf x})}{f^2({\bf X})}d{\bf x}.
\end{eqnarray}
But, using once again the arguments used to show
(\ref{argument2}), along with the facts that $G(\tau)>0$, $\psi$
and $g$ are bounded and $L$ is compactly supported, it holds that,
as $n\rightarrow\infty$, \vspace{-0.5cm}
\begin{eqnarray}\label{I2}
E(I^2)&=&\int E( {\epsilon_1^\star}^4|{\bf X}_{1}={\bf v}) \left(
\int L^2\left(\frac{{\bf x}-{\bf v}}{\ell_n}\right) \frac{g({\bf
x})}{f^2({\bf v})}d{\bf x}\right)^2f({\bf v}) d{\bf v},\nonumber\\
&=&\mathcal{O}(\ell_n^{2d})\big(1+\mathcal{O}(\rho_n)\big).
\end{eqnarray}
By (\ref{vartn1}) and (\ref{I2}), it follows that ${\rm Var} \
T_{n,1}^{1^\star} =\mathcal{O}(n^{-3}\ell_n^{-2d})$. Using the
Bienayme-Tchebytchev  inequality, we infer that, for all
$\varepsilon>0$, $P(n\ell_n^{d/2}|T_{n,1}^{1^\star}-
ET_{n,1}^{1^\star}| \geq \varepsilon)  \leq
\varepsilon^{-2}n^2\ell_n^d\times {\rm Var}
T_{n,1}^{1^\star}=\mathcal{O}((n\ell^d_n)^{-1})$. Thus,
$T_{n,1}^{1^\star}=
(n\ell_n^d)^{-1}B+o_p[({n\ell_n^{d/2}})^{-1}]$, which is
(\ref{1}).\vskip5pt \noindent{\bf Proof of (\ref{2}):}~ For $1\leq
i\leq n$, set ${\bf {\zeta}}_i=({\bf X}_i, \epsilon_i)$ and
$u_i({\bf x})=L[({\bf x}-{\bf X}_i)/\ell_n]$. Further introduce
${T_{n,2}^{1}}= \ell_n^{-2d}\!\!\int\! \sum_{i < j} u_i({\bf
x})u_j({\bf x})\epsilon_i\epsilon_j g({\bf x})f^{-1}({\bf
X}_i)f^{-1}({\bf X}_j) d{\bf x}$. Note that, in view of
(\ref{LILfoldes}) and (\ref{ango}), the dominated convergence
theorem ensures that ${T_{n,2}^{1^\star}}= {T_{n,2}^1}\!\!' +
\mathcal{O}(\rho_n)$ almost surely as $n\rightarrow\infty$, with
${T_{n,2}^1}\!\!':=2n^{-2}{T_{n,2}^{1}}$.
\\To establish (\ref{2}), we will make use of a central limit theorem
for U-statistics due to \cite{HALL84}. Set $M_n(\zeta_i,{\bf
\zeta}_j)=\ell_n^{-2d}\int u_i({\bf x})u_j({\bf
x}))\epsilon_i\epsilon_jg({\bf x}))f^{-1}({\bf X}_i)f^{-1}({\bf
X}_j)d{\bf x}$ and $N_n({\bf u}, {\bf v})=E(M_n({\bf \zeta}_1,{\bf
u}) M_n({\bf \zeta}_1,{\bf v}))$. To apply Hall's theorem to
${T_{n,2}^{1}}\!\!'$, the conditions $[T1], [T2]$ and $[T3]$ below
must be verified.
\begin{tabbing}
$[T_1]\;\;$\=\ $E\{M_n(\zeta_1,\zeta_2)|\zeta_1\}=0.$\\
$[T_2]$\>\ $E\{M_n^2(\zeta_1,\zeta_2)\}<\infty.$\\
$[T_3]$\>\ $|E\{N_n^2(\zeta_1,\zeta_2)\}+
n^{-1}E\{M_n^4(\zeta_1,\zeta_2)\}|\ /\
|E\{M_n^2(\zeta_1,\zeta_2)\}|^2
    \rightarrow 0$.
\end{tabbing}
$[T_1]$ is readily satisfied by making use of conditioning
arguments. Moreover, arguing as before,  the statement $[T_2]$
follows from routine analysis. To establish $[T_3]$, it is
sufficient to prove the results (\ref{Gn2}), (\ref{Hn2}) and
(\ref{hn42}) below. \vspace{-0.5cm}
\begin{eqnarray}
E(N_n^2(\zeta_1,\zeta_2) )& = & \mathcal{O}(1)\quad\mbox{ as } n\rightarrow\infty, \label{Gn2}\\
E(M_n^2(\zeta_1,
\zeta_2))& =& \ell_n^{-d}\frac{V}{2}+o(\ell_n^{-d})\quad\mbox{ as } n\rightarrow\infty,\label{Hn2}\\
 E(M_n^4(\zeta_1,
\zeta_2))& =&\mathcal{O}\big(\ell_n^{-3d}\big)\quad\mbox{ as }
n\rightarrow\infty.\label{hn42}
\end{eqnarray}
{\em Proof of (\ref{Gn2}):}~ Denote by $f_{{\bf X},\epsilon}$ the
joint density of $({\bf X}, \epsilon)$ (the existence of which
being ensured by the assumption (C.1), since $({\bf X},Y)$ is
supposed to have a joint density). It holds that \vspace{-0.5cm}
\begin{eqnarray*}
&&E\big[N_n^2(\zeta_1,\zeta_2)\big] = \frac{1}{\ell_n^{8d}}\int
... \int \Big[\int ... \int L\Big(\frac{{\bf x}_1-{\bf
v}_1}{\ell_n}\Big) L\Big(\frac{{\bf
x}_1-\eta_1}{\ell_n}\Big)L\Big(\frac{{\bf x}_2-{\bf
v}_1}{\ell_n}\Big)\\
&&\!\!\!\!\times  L\Big(\frac{{\bf
x}_2-\eta_2}{\ell_n}\Big)\frac{g({\bf x}_1)}{f^2({\bf
v}_1)}\frac{g({\bf x}_2)}{f^2({\bf v}_2)} f_{{\bf
X},\epsilon}({\bf v}_1,e_1)e_1^2e_2 e_3 d{\bf v}_1de_1d{\bf
x}_1d{\bf x}_2\Big]^2\!\!de_2de_3d\eta_1d\eta_2.
\end{eqnarray*}
Using classical changes of variables, together with the assumption $(G.1)$
and the fact that $L$ is compactly supported, (\ref{Gn2}) is
straightforward.\vskip5pt

\noindent{\em Proof of (\ref{Hn2}).}~ Set $\sigma_n^2 = \int\!\!
\int M_n^2(\omega_1, \omega_2)f_{{\bf
X},\epsilon}(\omega_1)f_{{\bf X},\epsilon}(\omega_2)
d\omega_1d\omega_2,$. Then, \vspace{-0.5cm}
\begin{eqnarray*}
\sigma_n^2 &= &\int \int \int \int \int \int \frac{1}{\ell_n^{4d}}
L\Big(\frac{{\bf x}_1-{\bf v}_1}{\ell_n}\Big)L\Big(\frac{{\bf
x}_1-{\bf v}_2}{\ell_n}\Big) L\Big(\frac{{\bf x}_2-{\bf
v}_1}{\ell_n}\Big)L\Big(\frac{{\bf x}_2-{\bf v}_2}{\ell_n}\Big)\\
&  & \times e_1^2 e_2^2 f_{{\bf X}_1,\epsilon_1} ({\bf v}_1,
e_1)f_{{\bf X}_2,\epsilon_2} ({\bf v}_2, e_2) \frac{g({\bf
x}_2)}{f^2({\bf v_2})}\frac{g({\bf x}_1)}{f^2({\bf v_1})}d{\bf
x}_1d{\bf x}_2 d{\bf v}_1d{\bf v}_2 de_1de_2.
\end{eqnarray*}
Next, noting that $\int e_1^2 f_{{\bf X},\epsilon} ({\bf v}_1,e_1)
de_1 = E(\epsilon_i^2|{\bf X}_i={\bf v}_1) f({\bf v}_1)$, it
follows that \vspace{-0.5cm}
\begin{eqnarray*}
&&\sigma_n^2 = \int \int \int \int  \frac{1}{\ell_n^{4d}}
L\Big(\frac{{\bf x}_1-{\bf v}_1}{\ell_n}\Big)L\Big(\frac{{\bf
x}_1-{\bf v}_2}{\ell_n}\Big)L\Big(\frac{{\bf x}_2-{\bf
v}_1}{\ell_n}\Big)L\Big(\frac{{\bf x}_2-{\bf v}_2}{\ell_n}\Big)\\
&&\!\!\!\!\!\!\times\ E(\epsilon_i^2|{\bf X}_i={\bf v}_1) f({\bf
v}_1)E(\epsilon_i^2|{\bf X}_j={\bf v}_2) f({\bf v}_2) \frac{g({\bf
x}_2)}{f^2({\bf v_2})}\frac{g({\bf x}_1)}{f^2({\bf v_1})}d{\bf
x}_1d{\bf x}_2 d{\bf v}_1d{\bf v}_2.
\end{eqnarray*}
Using  the changes of variables, ${\bf y}_1=({\bf x}_1-{\bf
v}_1)/\ell_n$, ${\bf y}_2=({\bf x}_2-{\bf v}_1)/\ell_n$ and ${\bf
r}_1=({\bf v}_1-{\bf v}_2)/\ell_n$, along with the continuity of
$f$ and the dominated convergence theorem, we get, \vspace{-0.5cm}
\begin{eqnarray*}
\sigma_n^2 &= &\int ... \int  \frac{1}{\ell_n^{4d}} L({\bf
y}_1)L({\bf y}_2)L({\bf y}_1-{\bf r}_1)L({\bf y}_2-{\bf
r}_1)E(\epsilon_i^2|{\bf X}_i={\bf v}_1)\\& &\times
E(\epsilon_j^2|{\bf X}_j={\bf v}_1+{\bf r}_1\ell_n)\frac{g({\bf
v}_1+{\bf y}_1\ell_n)}{f({\bf v}_1)}\frac{g({\bf v}_1+{\bf
y}_1\ell_n)}{f({\bf v}_1+{\bf
r}_1\ell_n)}d{\bf y}_1d{\bf y}_2 d{\bf r}_1d{\bf v}_1,\\
& = &\ell_n^{-d}\int \left[\int L({\bf t})L({\bf t}- {\bf r})
d{\bf t}\right]^2d{\bf r}\int (E(\epsilon_i^2|{\bf X}_i={\bf
r}))^2\frac{g^2({\bf r})}{f^2({\bf r})}d{\bf r} + o(\ell_n^{-d}),
\end{eqnarray*}
which, recalling the definition of $V$, is
(\ref{Hn2}).\vskip5pt\noindent {\em Proof of
(\ref{hn42})}:~Arguing as before (see also \cite{Camlong2001}), we
can show that, for a given $C<\infty$, \vspace{-0.5cm}
\begin{eqnarray*}
|E(M_n^4(\zeta_i,\zeta_j))|& \leq & \frac{C}{\ell_n^{8d}}\int
...\int\Big|L\Big(\frac{{\bf x}_1-{\bf
v}_1}{\ell_n}\Big)L\Big(\frac{{\bf x}_2-{\bf v}_1}{\ell_n}\Big)
\\&&L\Big(\frac{{\bf x}_3-{\bf v}_1}{\ell_n}\Big)
L\Big(\frac{{\bf x}_4-{\bf v}_1}{\ell_n}\Big)L\Big(\frac{{\bf
x}_1-{\bf v}_2}{\ell_n}\Big)L\Big(\frac{{\bf x}_2-{\bf
v}_2}{\ell_n}\Big)\\&&L\Big(\frac{{\bf x}_3-{\bf
v}_2}{\ell_n}\Big)L\Big(\frac{{\bf x}_4-{\bf
v}_2}{\ell_n}\Big)\Big|\frac{g({\bf x}_1)}{f({\bf
v}_1)}\frac{g({\bf x}_2)}{f({\bf v}_2)}\\&&\frac{g({\bf
x}_3)}{f^2({\bf v}_1)}\frac{g({\bf x}_4)}{f^2({\bf v}_2)} d{\bf
v}_1d{\bf v}_2d{\bf x}_1d{\bf x}_2d{\bf x}_3d{\bf x}_4
\\&=& \mathcal{O}\big(\ell_n^{-3d}\big).
\end{eqnarray*}
Combining (\ref{Gn2}), (\ref{Hn2}) and (\ref{hn42}), it is readily
shown that $[T3]$ holds. Then, Hall's Theorem can be applied to
${T_{n,2}^{1}}$. Namely, since $EM_n(\zeta_i, \zeta_j)=0$, we have
$\sqrt{2}{(n\sigma_n)^{-1}}T_{n,2}^{1'}\rightarrow
\mathcal{N}(0,1)$. Recalling that $T_{n,2}^{1}=n^2
T_{n,2}^{1'}/2$, we deduce, from (\ref{Hn2}), that
$(n\ell_n^{d/2})T_{n,2}^{1}/\sqrt{V}\rightarrow \mathcal{N}(0,1)$.
Slutsky's Theorem is now sufficient to conclude to (\ref{2}),
since, as already mentioned, ${T_{n,2}^{1^\star}}= T_{n,2}^1 +
\mathcal{O}(\rho_n)$ almost surely as $n\rightarrow\infty$
.\vskip5pt \noindent {\bf Proof of (\ref{3}):} By $(G.1)$,
\vspace{-0.5cm}
\begin{eqnarray*}
T_{n,3}^{1^\star}= \sup_{{\bf x}\in C}|m_\psi({\bf x})- \widehat
m_{\psi,add}^\star({\bf
x})|^2\int\Big[\frac{1}{n\ell_n^{d}}\sum_{i=1}^n L\Big(\frac{{\bf
x}-{\bf X}_i}{\ell_n}\Big)\Big]^2g({\bf x})d{\bf x}\quad {\rm
a.s.}\ .
\end{eqnarray*}
Since, under $H_0$, $m_\psi=m_{\psi,add} \in \mathcal{M}_{add}$,
we can apply the result of Lemma \ref{lem_unif_consist}.  This
latter, when combined with, successively, the boundedness of $f$,
the dominated convergence theorem, Bochner's theorem and the fact
that $g$ is compactly supported, yields $T_{n,3}^{1^\star} =
{\mathcal O}\big((\log n/n)^{2k/2k+1}\big)$ almost surely as
$n\rightarrow\infty$. Thus, under the assumption $(H.3)$,
$T_{n,3}^{1^\star}=o_{p}(n^{-1}\ell_n^{-d/2})$ as
$n\rightarrow\infty$.\vskip5pt

\noindent {\bf Proof of (\ref{4}):}~ First consider the mean of
$T_{n,4}^{1^\star}$. By $(G.1)$, it holds that \vspace{-0.5cm}
\begin{eqnarray*}
|ET_{n,4}^{1^\star}| 
& \leq & \sup_{{\bf x}\in C}|m_\psi({\bf x})- \widehat
m_{\psi,add}^\star({\bf x})|\frac{1}{n\ell_n^{2d}}\int \int
\Big|E\Big(\epsilon_i^\star|{\bf X}_i={\bf u}\Big)  L^2
\Big(\frac{{\bf x}-{\bf u}}{\ell_n} \Big)\Big|\\
&& \quad\times g({\bf x})d{\bf x}d{\bf u}.
\end{eqnarray*}
Next, under the assumptions $(H.1$-$2$-$3)$, using, successively,
the assumption $(G.1)$, the dominated convergence theorem, the
equality (\ref{argument2}), Bochner's theorem and Lemma
\ref{lem_unif_consist}, it can be shown that $|ET_{n,4}^{1^\star}|
= \mathcal{O}((n\ell_n^d)^{-1}(\log n/n)^{k/(2k+1)})$ $=
o(n^{-1}\ell_n^{-d/2})$. Turning our attention to the variance of
$T_{n,4}^{1^\star}$, and arguing as before, we get \vspace{-0.5cm}
\begin{eqnarray*}
{\rm Var}\ T_{n,4}^{1^\star}& = &\frac{1}{n^3\ell_n^{4d}}{\rm
Var}\Big[\int \epsilon_i^\star L^2\Big(\frac{{\bf x}-{\bf
X}_i}{\ell_n} \Big)(m_\psi({\bf X}_i)- \widehat
m_{\psi,add}^\star({\bf
X}_i))\frac{g({\bf x})}{f^2({\bf X}_i)}d{\bf x}\Big]\nonumber\\
& =&\mathcal{O}\Big(\frac{1}{n^3\ell_n^{2d}}\Big(\frac{\log
n}{n}\Big)^{2k/(2k+1)}\Big).
\end{eqnarray*}
An application of Bienayme-Tchebychev's inequality leads to
$T_{n,4}^{1^\star} - ET_{n,4}^{1^\star} =
o_p(n^{-1}\ell_n^{-d/2})$, which implies (\ref{4}), since
$|ET_{n,4}^{1^\star}|=o(n^{-1}\ell_n^{-d/2})$. \vskip5pt\noindent
{\bf Proof of (\ref{5}):} Arguing as before, we infer that,
ultimately as $n\rightarrow\infty$, \vspace{-0.5cm}
\begin{eqnarray*}
ET_{n,5}^{1^\star} = \mathcal{O}\Big(\Big(\frac{\log n}{n}\Big)^{k/(2k+1)}\Big)
\quad\mbox{and}\quad
{\rm Var}T_{n,5}^{1^\star}=\mathcal{O} \Big(\frac{1}{n^2}\Big(\frac{\log
n}{n}\Big)^{2k/(2k+1)}\Big).
\end{eqnarray*}
Therefore, Bienayme-Tchebychev's inequality
leads to (\ref{5}). 
%


\end{document}